\documentclass[preprint,1p]{elsarticle}
\usepackage{graphics}
\usepackage{epsfig}
\usepackage{mathptmx}
\usepackage{amsmath}
\usepackage{amssymb}
\usepackage{hyperref}
\usepackage{fullpage}
\usepackage{amsthm}
\usepackage{lineno}

\journal{Name of Journal}

\begin{document}
\begin{frontmatter}
\title{Approximation properties by generalized-Baskakov-Kantorovich-Stancu type operators}
\author[]{Abdul wafi}
\ead{abdulwafi2k2@gmail.com}
\author{Nadeem Rao\corref{cor1}}
\ead{nadeemrao1990@gmail.com}
\address{\textbf{Department of Mathematics, Jamia Millia Islamia, New Delhi-110025, India}}
\cortext[cor1]{Corresponding Author}
\begin{abstract}
  In this paper, we introduce generalized Baskakov Kantorovich Stancu type operators and investigate direct result, local approximation and weighted approximation properties of these operators. Modulus of continuity, second modulus of continuity, Peeter's K-functional, weighted modulus of continuity and Lipschitz class are considered to prove our results.
\end{abstract}
\end{frontmatter}
\section{Introduction}
    For $f\in C[0,1]$, Bernstein\cite{Bernstein} defined the linear positive operators which are the classical example of linear approximation as
\begin{eqnarray}
B_n(f;x)=\sum_{k=0}^{n}(^n_k)P_{n,k}(x)f\bigg(\frac{k}{n}\bigg),
\end{eqnarray}
where $P_{n,k}(x)=x^k(1-x)^{n-k}, x\in[0,1]$ and $n\in N$. In this paper, Bernstein showed that these operators approximate uniformly on [0,1] to every continuous function $f\in C[0,1].$ But these operators are not suitable for discontinuous functions. Later on, Kantorovich\cite{kan} generalized the operators(1) to approximate the measurable functions (see G.G. Lorenz\cite{lorentz}). For $n\in N$ and $f\in L_p[0,1], 1\leq p<\infty$, the Kantorovich operators $K_n:L_p([0,1])\rightarrow L_p([0,1])$ defied by
\begin{eqnarray}
  K_n(f;x)=(n+1)\sum_{k=0}^{n}P_{n,k}(x)\int\limits_{\frac{k}{n+1}}^{\frac{k+1}{n+1}}f(t)dt, \hspace{2cm} k=1,2,3...
\end{eqnarray}
Kantorovich type operators converges almost everywhere to the function $f$ on [0,1]. A. Wafi and S. Khatoon\cite{wafi2} defined a generalized Baskakov-Kantorovich operators on $[0,\infty)$
\begin{eqnarray}
V_n^a(f;x)=n \sum_{k=0}^{\infty}W_{n,k}^a(x)\int \limits_{\frac{k}{n}}^{\frac{k+1}{n}}f(t)dt,
\end{eqnarray}
where
\begin{eqnarray}
W_{n,k}^a(x)=e^{-\frac{ax}{1+x}}\frac{p_k(n,a)}{k!} \frac{x^k}{(1+x)^{k+n}}
\end{eqnarray}
such that $\sum\limits_{k=0}^{\infty}W_{n,k}^a(x)=1$ and $p_k(n,a)=\sum\limits_{i=0}^{\infty}\big( ^n _i\big)(n)_ia^{k-i},$ with $(n)_0=1,
 (n)_i=n(n+1)...(n+i-1).$\\ \\
The Stancu operators\cite{D.D.} and its Kantorovich version \cite{barbosu} are respectively given by
\begin{eqnarray}
B_n^{\alpha,\beta}(f;x)=(n+\beta+1)\sum_{k=0}^{n}x^k(1-x)^{n-k}f\bigg(\frac{k+\alpha}{n+\beta}\bigg)\\
T_n^{\alpha,\beta}(f;x)=(n+\beta+1)\sum_{k=0}^{n}x^k(1-x)^{n-k}\int\limits_{\frac{k+\alpha}{n+\beta+1}}^{\frac{k+\alpha+1}{n+\beta+1}}f(t)dt,
\end{eqnarray}
where $\alpha,\beta$ are any two two non negative real numbers such that $0\leq\alpha\leq\beta$. If $\alpha=\beta=0$, the operators (5) and (6) reduce to operators (1) and (2) respectively. In the literature there are many studies on these type of operators(see \cite{aral}, \cite{barbosu}, \cite{guptab}, \cite{guptads}, \cite{lopez}, \cite{sucu}). In this article, we define generalized Baskakov-Kantorovich-Stancu type opeators as follows
\begin{eqnarray}
T_{n,a}^{\alpha,\beta}(f;x)=(n+\beta) \sum_{k=0}^{\infty}W_{n,k}^a(x)\int \limits_{\frac{k+\alpha}{n+\beta}}^{\frac{k+\alpha+1}{n+\beta}}f(t)dt,
\end{eqnarray}
where $W_{n,k}^a(x)$ is defined in (4).\\ \\
For $\alpha=\beta=0$, we get the operators (3). The aim of this article is to present a Kantorovich form of Stancu variant of generaliged Baskakov operators (see Rao. N., Wafi A.\cite{nadeem}). Further, we study the direct result, local approximation weighted Korovkin type theorem and order of approximation using weighted modulus of continuity.
\section{Basic results }
\textbf{Lemma 2.1} Let $a,x\geq0$ and $n=1,2,3....$ Then from\cite{nadeem}, we have
\begin{eqnarray*}
(i)\quad L_{n,a}^{\alpha,\beta}(1;x)&=&1,\\
(ii) \quad L_{n,a}^{\alpha,\beta}(t;x)&=&\frac{n}{n+\beta}x+\frac{a}{n+\beta}\frac{x}{1+x}+ \frac{\alpha}{n+\beta},\\
(iii)\quad L_{n,a}^{\alpha,\beta}(t^2;x)&=&\frac{n^2+n}{(n+\beta)^2}x^2+\frac{n(1+2\alpha)}{(n+\beta)^2}x+\frac{a^2}{(n+\beta)^2}\frac{x^2}{(1+x)^2}+\frac{2an}{(n+\beta)^2}\frac{x^2}{(1+x)}+\frac{a(1+2\alpha)}{(n+\beta)^2}\frac{x}{1+x} + \frac{\alpha^2}{(n+\beta)^2}\\
(iv) \quad L_{n,a}^{\alpha,\beta}(t^3;x)&=&\frac{n^3+3n^2+2n}{(n+\beta)^3}x^3+\frac{n^2(3+3\alpha)+n(3+3\alpha+3a)}{(n+\beta)^3}x^2+\frac{n(1+3\alpha+3\alpha^2)}{(n+\beta)^3}x+\frac{3an^2}{(n+\beta)^3}\frac{x^3}{(1+x)}\\
&&+\frac{n}{(n+\beta)^3}\bigg(\frac{3a^2x^3}{(1+x)^2}+\frac{3ax^2}{1+x}+\frac{6a \alpha x^2}{1+x} \bigg)+\frac{1}{(n+\beta)^3}\bigg(\frac{ax}{1+x}+\frac{3a^2x^2}{(1+x)^2}+\frac{a^3x^3}{(1+x)^3}\\
&&+\frac{3\alpha a^2x^2}{(1+x)^2}+\frac{3\alpha^2 ax}{1+x} +\alpha^3\bigg),\\
(v) \quad L_{n,a}^{\alpha,\beta}(t^4;x)&=&\frac{n^4+6n^3+11n^2+6n}{(n+\beta)^4}x^4+\frac{(6+4\alpha)n^3+(18+12\alpha)n^2+(9+8\alpha)n}{(n+\beta)^4}x^3+\bigg(\frac{(7+12\alpha+6\alpha^2)n^2}{(n+\beta)^4}\\
&&+\frac{(7+12\alpha+12\alpha a+6\alpha ^2)n}{(n+\beta)^4}\bigg)x^2+\frac{(1+4\alpha+6\alpha ^2+4\alpha ^3)n}{(n+\beta)^4}x+\frac{4an^3+12an^2+8an}{(n+\beta)^4}\frac{x^4}{1+x}\\
&&+\frac{6a^2n^2+6a^2n}{(n+\beta)^4}\frac{x^4}{(1+x)^2}+\frac{4a^3n}{(n+\beta)^4}\frac{x^4}{(1+x)^3}+\frac{a^4}{(b+\beta)^4}\frac{x^4}{(1+x)^4}+\frac{18an^2+18an}{(n+\beta)^4}\frac{x^3}{(1+x)}\\
&&+\frac{(18a^2+12a^2\alpha)n}{(n+\beta)^4}\frac{x^3}{(1+x)^2}+\frac{6a^3+4\alpha a^3}{(n+\beta)^4}\frac{x^3}{(1+x)^3}+\frac{(12a\alpha ^2+12a\alpha+14a)n}{(n+\beta)^4}\frac{x^2}{(1+x)}\\
&&+\frac{7a^2+12a^2\alpha+6a^2\alpha^2}{(n+\beta)^4}\frac{x^2}{(1+x)^2}+\frac{a+4\alpha a+6\alpha^2a+4\alpha^3a}{(n+\beta)^4}\frac{x}{1+x}+\frac{\alpha^4}{(n+\beta)^4}.
\end{eqnarray*}
\textbf{Lemma 2.2} Let $a,x\geq0$. Then
\begin{eqnarray*}
(i)\quad T_{n,a}^{\alpha,\beta}(1;x)&=&1,\\
(ii) \quad T_{n,a}^{\alpha,\beta}(t;x)&=&\frac{n}{n+\beta}x+\frac{a}{n+\beta}\frac{x}{1+x}+ \frac{2\alpha+1}{2(n+\beta)},\\
(iii)\quad T_{n,a}^{\alpha,\beta}(t^2;x)&=&\frac{n^2+n}{(n+\beta)^2}x^2+\frac{n(2+2\alpha)}{(n+\beta)^2}x+\frac{a^2}{(n+\beta)^2}\frac{x^2}{(1+x)^2}+\frac{2an}{(n+\beta)^2}\frac{x^2}{(1+x)}+\frac{a(2+2\alpha)}{(n+\beta)^2}\frac{x}{1+x} + \frac{3\alpha^2+1}{3(n+\beta)^2}\\
(iv) \quad T_{n,a}^{\alpha,\beta}(t^3;x)&=&\frac{n^3+3n^2+2n}{(n+\beta)^3}x^3+\frac{n^2(\frac{9}{2}+3\alpha)+n(\frac{9}{2}+3\alpha+3a)}{(n+\beta)^3}x^2+\frac{n(\frac{7}{2}+6\alpha+3\alpha^2)}{(n+\beta)^3}x+\frac{3an^2}{(n+\beta)^3}\frac{x^3}{(1+x)}\\
&&+\frac{n}{(n+\beta)^3}\bigg(\frac{3a^2x^3}{(1+x)^2}+\frac{6ax^2}{1+x}+\frac{6a \alpha x^2}{1+x} \bigg)+\frac{1}{(n+\beta)^3}\bigg(a\big(\frac{7}{2}+3\alpha+3\alpha a^2\big)\frac{x}{1+x}+\big(\frac{9}{2}a^2+3\alpha a^2\big)\frac{x^2}{(1+x)^2}\\
&&+\frac{a^3x^3}{(1+x)^3}\bigg)+\frac{\alpha^3+\frac{3}{2}\alpha ^2+\alpha}{(n+\beta)^3},\\
(v) \quad T_{n,a}^{\alpha,\beta}(t^4;x)&=&\frac{n^4+6n^3+11n^2+6n}{(n+\beta)^4}x^4+\frac{(8+4\alpha)n^3+(24+12\alpha)n^2+(13+8\alpha)n}{(n+\beta)^4}x^3+\bigg(\frac{(15+18\alpha+6\alpha^2)n^2}{(n+\beta)^4}\\
&&+\frac{(15+18\alpha++6a+12\alpha a+6\alpha ^2)n}{(n+\beta)^4}\bigg)x^2+\frac{(6+14\alpha+11\alpha ^2+4\alpha ^3)n}{(n+\beta)^4}x+\frac{4an^3+12an^2+8an}{(n+\beta)^4}\frac{x^4}{1+x}\\
&&+\frac{6a^2n^2+6a^2n}{(n+\beta)^4}\frac{x^4}{(+x)^2}+\frac{4a^3n}{(n+\beta)^4}\frac{x^4}{(1+x)^3}+\frac{a^4}{(b+\beta)^4}\frac{x^4}{(1+x)^4}+\frac{24an^2+18an}{(n+\beta)^4}\frac{x^3}{(1+x)}\\
&&+\frac{(24a^2+12a^2\alpha)n}{(n+\beta)^4}\frac{x^3}{(1+x)^2}+\frac{8a^3+4\alpha a^3}{(n+\beta)^4}\frac{x^3}{(1+x)^3}+\frac{(12a\alpha ^2+10a+24a\alpha+14a)n}{(n+\beta)^4}\frac{x^2}{(1+x)}\\
&&+\frac{15a^2+18a^2\alpha+6a^2\alpha^2}{(n+\beta)^4}\frac{x^2}{(1+x)^2}+\frac{6a+8\alpha a+12\alpha^2a+4\alpha^3a}{(n+\beta)^4}\frac{x}{1+x}+\frac{\alpha^4+2\alpha^3+2\alpha^2+\alpha}{(n+\beta)^4}.
\end{eqnarray*}
\textbf{Proof} Identity (i) is obvious. For identity (ii), we use the linearity property
\begin{eqnarray*}
T_{n,a}^{\alpha,\beta}(t;x)=(n+\beta) \sum_{k=0}^{\infty}W_{n,k}^a(x)\int \limits_{\frac{k+\alpha}{n+\beta}}^{\frac{k+\alpha+1}{n+\beta}}tdt,
\end{eqnarray*}
\begin{eqnarray*}
T_{n,a}^{\alpha,\beta}(f;x)=(n+\beta) \sum_{k=0}^{\infty}W_{n,k}^a(x)\bigg[\frac{t^2}{2}\bigg]_{\frac{k+\alpha}{n+\beta}}^{\frac{k+\alpha+1}{n+\beta}},
\end{eqnarray*}
\begin{eqnarray*}
T_{n,a}^{\alpha,\beta}(t;x)=\sum_{k=0}^{\infty}W_{n,k}^a(x)\bigg[\frac{k+\alpha}{n+\beta}+\frac{1}{2(n+\beta)}\bigg].
\end{eqnarray*}
From Lemma 2.1, we have
\begin{eqnarray*}
 T_{n,a}^{\alpha,\beta}(t;x)=L_{n,a}^{\alpha,\beta}(t;x)+\frac{1}{2(n+\beta)}L_{n,a}^{\alpha,\beta}(1;x)
\end{eqnarray*}
\begin{eqnarray*}
T_{n,a}^{\alpha,\beta}(t;x)=\frac{n}{n+\beta}x+\frac{a}{n+\beta}\frac{x}{1+x}+ \frac{\alpha}{n+\beta}+\frac{1}{2(n+\beta)}
\end{eqnarray*}
which prove the identity (ii). Similarly, we can prove identities (iii), (iv) and (v).\\ \\
\textbf{Lemma 2.3} Let $\psi_x^i(t)=(t-x)^i, i=1,2,3,...$. Then, we have
\begin{eqnarray*}
T_{n,a}^{\alpha,\beta}(\psi_x^0(t);x)&=&1,\\
T_{n,a}^{\alpha,\beta}(\psi_x^1(t);x)&=&\bigg(\frac{n}{n+\beta}-1\bigg)x+\frac{a}{n+\beta}\frac{x}{1+x}+ \frac{2\alpha+1}{2(n+\beta)},\\
T_{n,a}^{\alpha,\beta}(\psi_x^2(t);x)&=&\frac{n+\beta^2}{(n+\beta)^2}x^2+\frac{n-2(\alpha+1)\beta}{(n+\beta)^2}x+\frac{a^2}{(n+\beta)^2}\frac{x^2}{(1+x)^2}-\frac{2a\beta}{(n+\beta)^2}\frac{x^2}{(1+x)}+\frac{a(2+2\alpha)}{(n+\beta)^2}\frac{x}{1+x} + \frac{3\alpha^2+1}{3(n+\beta)^2},\\
T_{n,a}^{\alpha,\beta}(\psi_x^4(t);x)&=&\frac{(3-12\beta)n^2+(6+4\beta+2\beta^2+4\beta^3)n+\beta^4}{(n+\beta)^4}x^4+\bigg(\frac{(6-12a-12\beta\alpha)n^2+(13+8\alpha-18\beta+12a\beta+9\beta^2)n}{(n+\beta)^4}\\
&&+\frac{-4\beta^3\alpha-2\beta^2}{(n+\beta)^4}\bigg)x^3+\frac{(11-18\alpha+18\alpha^2)n^2+(15+18\alpha+6a+12\alpha a-24 \alpha \beta)n+6\alpha^2\beta^2+2\beta^2}{(n+\beta)^4}x^2\\
&&+\frac{(6+10\alpha+5\alpha^2)n-4\beta(\alpha^3+\frac{3}{2}\alpha^2+\alpha)}{(n+\beta)^4}x+\frac{12an^2+8an-4a\beta^3}{(n+\beta)^4}\frac{x^4}{(1+x)}\\
&&+\frac{6a^2n+6a^2\beta^2}{(n+\beta)^4}\frac{x^4}{(1+x)^2} -\frac{4a^3\beta}{(n+\beta)^4}\frac{x^4}{(1+x)^3}+\frac{a^4}{(n+\beta)^4}\frac{x^4}{(1+x)^4}+\frac{12an^2+18an+6a(1+2\alpha)\beta^2}{(n+\beta)^4}\frac{x^3}{1+x}\\
&&+\frac{6a^2n-(12a^2+12\alpha a^2)\beta}{(n+\beta)^4}\frac{x^3}{(1+x)^2}+\frac{6a^3+4\alpha a^3}{(n+\beta)^4}\frac{x^3}{(1+x)^3}+\frac{(12a\alpha+8a-6a\alpha^2)n-(6a+18\alpha^2a)\beta}{(n+\beta)^4}\frac{x^2}{1+x}\\
&&+\frac{7a^2+12a^2\alpha+6a^2\alpha^2}{(n+\beta)^4}\frac{x^2}{(1+x)^2}+\frac{(a)+4\alpha a+6\alpha^2a+4\alpha^3a}{(n+\beta)^4}\frac{x}{1+x}+\frac{\alpha^4}{(n+\beta)^4}. \\
\end{eqnarray*}
\textbf{Proof} By using the linearity property and Lemma 2.2, we can prove the Lemma2.3 .\\ \\
\textbf{Lemma 2.4} For the operators $T_{n,a}^{\alpha,\beta}$, we have
\begin{eqnarray*}
 T_{n,a}^{\alpha,\beta}((t-x)^4;x)&\leq& M_a^{\alpha,\beta}\frac{(x^4+x^3+x^2+x+1)}{(n+\beta)^2}
 \end{eqnarray*}
 where $M_a^{\alpha,\beta}$ is a positive constant for the fixed value of non-negative numbers $a,\alpha$, and $\beta$.\\ \\
 \textbf{Proof} Since $\frac{x^i}{(1+x)^j}\leq x^i$ for all $x\geq0$ and $j\leq i(j,i=1,2,3,4)$, we have
  \begin{eqnarray*}
T_{n,a}^{\alpha,\beta}((t-x)^4;x)&\leq& \frac{1}{(n+\beta)^2} \Bigg(\frac{(3-12\beta+12a)n^2+(6+4\beta+2\beta^2+4\beta^3+8a+6a^2)n+\beta^4-4a\beta^3+6a^2\beta^2-4a^3\beta+a^4}{(n+\beta)^2}x^4\\
 &&+\bigg(\frac{(6-12a-12\beta\alpha+12a)n^2+(13+8\alpha-18\beta+12a\beta+9\beta^2+18a+6a^2)n-4\beta^3\alpha-2\beta^2}{(n+\beta)^2}\\
&&\frac{ (6a(1+2\alpha)\beta^2-12a^2+12\alpha a^2)\beta+6a^3+4\alpha a^3}{(n+\beta)^4}\bigg)x^3+\bigg(\frac{(11-18\alpha+18\alpha^2)n^2}{(n+\beta)^2}\\
&&+\frac{(15+18\alpha+14a+24\alpha a-24\alpha\beta-6a\alpha^2)n+12\alpha^2\beta^2+2\beta^2-(6a+18\alpha^2a)\beta+7a^2+12a^2\alpha}{(n+\beta)^2}\bigg)x^2\\
&&+\frac{(6+10\alpha+5\alpha^2)n-4\beta(\alpha^3+\frac{3}{2}\alpha^2+\alpha)+a+4\alpha a+6\alpha^2a+6\alpha^3a}{(n+\beta)^2}x,
   \end{eqnarray*}
\begin{eqnarray}
T_{n,a}^{\alpha,\beta}((t-x)^4;x)\leq\frac{A_4 x^4+A_3 x^3+A_2 x^2+A_1 x }{(n+\beta)^2}
\end{eqnarray}
   where
   \begin{eqnarray*}
   A_4&=&\frac{(3-12\beta+12a)n^2+(6+4\beta+2\beta^2+4\beta^3+8a+6a^2)n+\beta^4-4a\beta^3+6a^2\beta^2-4a^3\beta+a^4}{(n+\beta)^2}\\
   A_3&=&\frac{(6-12a-12\beta\alpha+12a)n^2+(13+8\alpha-18\beta+12a\beta+9\beta^2+18a+6a^2)n-4\beta^3\alpha-2\beta^2}{(n+\beta)^2}\\
&&+\frac{ (6a(1+2\alpha)\beta^2-12a^2+12\alpha a^2)\beta+6a^3+4\alpha a^3}{(n+\beta)^4}\\
A_2&=&\frac{(15+18\alpha+14a+24\alpha A-24\alpha\beta-6a\alpha^2)n+12\alpha^2\beta^2+2\beta^2-(6a+18\alpha^2a)\beta+7a^2+12a^2\alpha}{(n+\beta)^2}\bigg)\\
A_1&=&\frac{(6+10\alpha+5\alpha^2)n-4\beta(\alpha^3+\frac{3}{2}\alpha^2+\alpha)+a+4\alpha a+6\alpha^2a+6\alpha^3a}{(n+\beta)^2}.
   \end{eqnarray*}
   For large value of $n$, each $A_i(i=1,2,3,4)$ converges to a value depending on $a,\alpha,\beta$, Therefore, there exists a constant $M_a^ {\alpha,\beta}$ for the fixed values of $a,\alpha,\beta$ such that $A_i\leq M_a^ {\alpha,\beta}$ for $i=1,2,3,4$. Using (8), we arrive at the result. \\

\section{\textbf{Direct approximation}}
Let $C_B[0,\infty)$ denote the space of real valued continuous and bounded functions $f$ on $[0,\infty)$ endowed with the norm
\begin{eqnarray*}
\|f\|=\sup\limits_{0\leq x<\infty}|f(x)|.
\end{eqnarray*}
Then, for any $\delta>0$, Peeter's K-functional is defined as
\begin{eqnarray*}
K_2(f,\delta)=inf\{\|f-g\|+\delta\|g''\|: g\in C_B^2[0,\infty)\},
\end{eqnarray*}
where $C_B^2[0,\infty)=\{g\in C_B[0,\infty):g',g''\in C_B[0,\infty)\}$. By Devore and Lorentz[\cite{devor}, p.177, Theorem 2.4], there exits an absolute constant $C>0$ such that
\begin{eqnarray*}
K_2(f;\delta)\leq C\omega_2(f;\sqrt{\delta}),
\end{eqnarray*}
where $\omega_2(f;\delta)$ is the second order modulus of continuity is defined as
\begin{eqnarray*}
\omega_2(f,\sqrt{\delta})=\sup\limits_{0<h\leq\sqrt{\delta}} \sup\limits_{x\in[0,\infty)} |f(x+2h)-2f(x+h)+f(x)|.
\end{eqnarray*}
\textbf{Theorem 3.1} Let $f\in C_B^2[0,\infty)$. Then for all $x\in[0,\infty)$ there exist a constant $K>0$ such that
\begin{eqnarray*}
\mid T_{n,a}^{\alpha,\beta}(f;x)-f(x)\mid\leq K\omega_2(f;\sqrt{\gamma_n^{\alpha,\beta}(x)})+\omega\Bigg(f;\bigg(\frac{\beta}{n+\beta}\bigg)x+\frac{a}{n+\beta}\frac{x}{1+x}+ \frac{2\alpha+1}{2(n+\beta)}\Bigg)
\end{eqnarray*}
where $\gamma_n^{\alpha,\beta}(x)=\Bigg\{\frac{n+2\beta^2}{(n+\beta)^2}x^2+\frac{n-\beta}{(n+\beta)^2}x+\frac{2a^2}{(n+\beta)^2}\frac{x^2}{(1+x)^2}+\frac{a(3+4\alpha)}{(n+\beta)^2}\frac{x}{1+x} + \frac{7\alpha^2+4\alpha+2}{3(n+\beta)^2}\Bigg\}$.\\ \\
\textbf{Proof} First, we define the auxiliary operators
\begin{eqnarray}
\hat{T}_n^{\alpha,\beta}(f;x)=T_{n,a}^{\alpha,\beta}(f;x)+f(x)-f\Big(\frac{n}{n+\beta}x+\frac{a}{n+\beta}\frac{x}{1+x}+ \frac{2\alpha+1}{2(n+\beta)}\Big)
\end{eqnarray}
We find that
\begin{eqnarray*}
\hat{T}_n^{\alpha,\beta}(1;x)=1,
\end{eqnarray*}
\begin{eqnarray*}
\hat{T}_n^{\alpha,\beta}(\psi_x(t);x)=0
\end{eqnarray*}
\begin{eqnarray}
|\Hat{T}_n^{\alpha,\beta}(f;x)|\leq 3\|f\|.
\end{eqnarray}
Let $g\in C_B^2[0,\infty)$. By the Taylor's theorem\\
\begin{eqnarray}
g(t)=g(x)+(t-x)g'(x)+\int\limits_x^t (t-v)g''(v)dv,
\end{eqnarray}
Now, using (11), the auxiliary operators is given
\begin{eqnarray*}
\Hat{T}_n^{\alpha,\beta}(g;x)-g(x)&=&g'(x)\Hat{T}_n^{\alpha,\beta}(t-x;x)+\Hat{T}_n^{\alpha,\beta}\Big( \int\limits_x^t (t-v)g''(v)dv;x\Big)\\
&=&\Hat{T}_n^{\alpha,\beta}\Big( \int\limits_x^t (t-v)g''(v)dv;x\Big)\\
&=&T_{n,a}^{\alpha,\beta}\Big( \int\limits_x^t (t-v)g''(v)dv;x\Big)-\int\limits_x^{\frac{n}{n+\beta}x+\frac{a}{n+\beta}\frac{x}{1+x}+ \frac{2\alpha+1}{2(n+\beta)}} \Big(\frac{n}{n+\beta}x+\frac{a}{n+\beta}\frac{x}{1+x}+ \frac{2\alpha+1}{2(n+\beta)}-v\Big)g''(v)dv.
\end{eqnarray*}
Therefore
\begin{eqnarray}
\nonumber | \Hat{T}_n^{\alpha,\beta}(g;x)-g(x)|&\leq&\Bigg|T_{n,a}^{\alpha,\beta}\Big( \int\limits_x^t (t-v)g''(v)dv;x\Big)\Bigg|\\
&&+\Bigg|\int\limits_x^{\frac{n}{n+\beta}x+\frac{a}{n+\beta}\frac{x}{1+x}+ \frac{2\alpha+1}{2(n+\beta)}} \Big(\frac{n}{n+\beta}x+\frac{a}{n+\beta}\frac{x}{1+x}+ \frac{2\alpha+1}{2(n+\beta)}-v\Big)g''(v)dv\Bigg|.
\end{eqnarray}
Since
 \begin{eqnarray}
\Bigg| \int\limits_x^t (t-v)g''(v)dv\Bigg|\leq(t-x)^2\parallel g''\parallel
 \end{eqnarray}
 and
 \begin{eqnarray}
 \Bigg|\int\limits_x^{\frac{n}{n+\beta}x+\frac{a}{n+\beta}\frac{x}{1+x}+ \frac{2\alpha+1}{2(n+\beta)}} \Big(\frac{n}{n+\beta}x+\frac{a}{n+\beta}\frac{x}{1+x}+ \frac{2\alpha+1}{2(n+\beta)}-v\Big)g''(v)dv\Bigg|\leq \Bigg(\frac{\beta}{n+\beta}x+\frac{a}{n+\beta}\frac{x}{1+x}+ \frac{2\alpha+1}{2(n+\beta)}\Bigg)^2\parallel g'' \parallel.
 \end{eqnarray}
 Then from(12),(13)and (14) implies that
 \begin{eqnarray}
\nonumber | \Hat{T}_n^{\alpha,\beta}(g;x)-g(x)|&\leq& \Bigg\{T_{n,a}^{\alpha,\beta}((t-x)^2;x)+\Bigg(\frac{\beta}{n+\beta}x+\frac{a}{n+\beta}\frac{x}{1+x}+ \frac{2\alpha+1}{2(n+\beta)}\Bigg)^2\Bigg\}\|g''\|\\
 \nonumber &=&\Bigg\{\frac{n+2\beta^2}{(n+\beta)^2}x^2+\frac{n-\beta}{(n+\beta)^2}x+\frac{2a^2}{(n+\beta)^2}\frac{x^2}{(1+x)^2}+\frac{a(3+4\alpha)}{(n+\beta)^2}\frac{x}{1+x} + \frac{7\alpha^2+4\alpha+2}{3(n+\beta)^2}\Bigg\}\|g''\|\\
 &=&\gamma_n^{\alpha,\beta}(x)\|g''\|
 \end{eqnarray}
 Next, we have
 \begin{eqnarray*}
 |T_{n,a}^{\alpha,\beta}(f;x)-f(x)|\leq |\Hat{T}_n^{\alpha,\beta}(f-g;x)|+|(f-g)(x)|+|\Hat{T}_n^{\alpha,\beta}(g;x)-g(x)|+\big|f(\frac{n}{n+\beta}x+\frac{a}{n+\beta}\frac{x}{1+x}+ \frac{2\alpha+1}{2(n+\beta)})-f(x)\big|
 \end{eqnarray*}
 Using(15), we have
 \begin{eqnarray*}
 |T_{n,a}^{\alpha,\beta}(f;x)-f(x)|&\leq& 4\|f-g\| +\Hat{T}_n^{\alpha,\beta}(g;x)-g(x)|+\big|f(\frac{n}{n+\beta}x+\frac{a}{n+\beta}\frac{x}{1+x}+ \frac{2\alpha+1}{2(n+\beta)})-f(x)\big|\\
 &\leq& 4\|f-g\|+\gamma_n^{\alpha,\beta}(x)\|g''\|+\omega\Big(f;\frac{n}{n+\beta}x+\frac{a}{n+\beta}\frac{x}{1+x}+ \frac{2\alpha+1}{2(n+\beta)}\Big)
 \end{eqnarray*}
 By the definition of Peetre's K-functional
 \begin{eqnarray*}
 |T_{n,a}^{\alpha,\beta}(f;x)-f(x)|\leq C\omega_2\big(f;\sqrt{\gamma_n^{\alpha,\beta}(x)}\big)+\omega\Big(f;\frac{n}{n+\beta}x+\frac{a}{n+\beta}\frac{x}{1+x}+ \frac{2\alpha+1}{2(n+\beta)}\Big).
 \end{eqnarray*}
 Now, we use the Lipschitz type space
 \begin{eqnarray*}
 LIp_M^*=\{f\in C[0,\infty):|f(t)-f(x)|\leq M\frac{|t-x|^\alpha}{(t+x)^{\frac{\alpha}{2}}}:x,t\in(0,\infty)\}
 \end{eqnarray*}
 where $M$ is a constant and $0<\alpha\leq 1$ to prove the following theorem:\\ \\ \\
 \textbf{Theorem 3.2} Let $f\in Lip_M^*(\alpha)$ and $x\in [0,\infty)$. Then, we have
 \begin{eqnarray*}
  |T_{n,a}^{\alpha,\beta}(f;x)-f(x)|\leq M\Bigg[\frac{\Lambda_n^{\alpha,\beta}(x)}{x}\Bigg],
 \end{eqnarray*}
 where $\Lambda_n^{\alpha,\beta}(x)=T_n^{\alpha,\beta}((t-x)^2;x)$. \\ \\
 \textbf{Proof} Let $\alpha=1$ and $x\in(0,\infty)$. Then, for $f\in Lip_M^*(1)$, we have
 \begin{eqnarray*}
  |T_{n,a}^{\alpha,\beta}(f;x)-f(x)|&\leq & (n+\beta)\sum_{k=0}^{\infty}W_{n,k}^a(x)\int_{\frac{k+\alpha}{n+\beta}}^{\frac{k+\alpha+1}{n+\beta}}|f(t)-f(x)|dt\\
  &\leq & M(n+\beta)\sum_{k=0}^{\infty}W_{n,k}^a(x)\int_{\frac{k+\alpha}{n+\beta}}^{\frac{k+\alpha+1}{n+\beta}}\frac{|t-x|}{\sqrt{t+x}}dt.\\
  &\leq & \frac{M}{\sqrt{x}}(n+\beta)\sum_{k=0}^{\infty}W_{n,k}^a(x)\int_{\frac{k+\alpha}{n+\beta}}^{\frac{k+\alpha+1}{n+\beta}}|t-x|dt\\
  &\leq&\frac{M}{\sqrt{x}}T_n^{\alpha,\beta}(|t-x|;x)\\
  &\leq& M \frac{\sqrt{T_n^{\alpha,\beta}((t-x)^2;x)}}{\sqrt{x}}\\
  &=&M \Bigg(\frac{\Lambda_n^{\alpha,\beta}(x)}{x}\Bigg)^{\frac{1}{2}}
 \end{eqnarray*}
 Thus, the assertion hold for $\alpha=1$. Now, we will prove for $\alpha\in (0,1)$. From the Holder inequality with $p=\frac{1}{\alpha}, q=\frac{1}{1-\alpha}$, we have
 \begin{eqnarray*}
  |T_{n,a}^{\alpha,\beta}(f;x)-f(x)|&=& \Bigg(\sum_{k=0}^{\infty}W_{n,k}^a(x)\Bigg((n+\beta)\int_{\frac{k+\alpha}{n+\beta}}^{\frac{k+\alpha+1}{n+\beta}}|f(t)-f(x)|dt\Bigg)^{\frac{1}{\alpha}}\Bigg)^\alpha\big(W_{n,k}^a(x)\big)^{1-\alpha}\\
  &\leq& \Bigg(\sum_{k=0}^{\infty}W_{n,k}^a(x)\Bigg((n+\beta)\int_{\frac{k+\alpha}{n+\beta}}^{\frac{k+\alpha+1}{n+\beta}}|f(t)-f(x)|dt\Bigg)^{\frac{1}{\alpha}}\Bigg)^\alpha
 \end{eqnarray*}
 Since $f\in Lip_M^*$, we obtain
 \begin{eqnarray*}
  |T_{n,a}^{\alpha,\beta}(f;x)-f(x)|&\leq& M\Bigg((n+\beta)\sum_{k=0}^{\infty}W_{n,k}^a(x)(\int_{\frac{k+\alpha}{n+\beta}}^{\frac{k+\alpha+1}{n+\beta}}\frac{|t-x|}{\sqrt{t+x}}dt\Bigg)^\alpha\\ &\leq&\frac{M}{x^{\frac{\alpha}{2}}}\Bigg((n+\beta)\sum_{k=0}^{\infty}W_{n,k}^a(x)\int_{\frac{k+\alpha}{n+\beta}}^{\frac{k+\alpha+1}{n+\beta}}|t-x|dt\Bigg)^\alpha\\
  &=&\frac{M}{x^{\frac{\alpha}{2}}}\big(T_{n,a}^{\alpha,\beta}(|t-x|;x)\big)^{\alpha}\\
  &\leq& M \Bigg(\frac{\Lambda_n^{\alpha,\beta}(x)}{x}\Bigg)^{\frac{r}{2}}.
  \end{eqnarray*}

\section{\textbf{Approximation properties of $T_{n,a}^{\alpha,\beta}$ in weighted space}}
In this section, we investigate the weighted approximation properties of the operators $T_{n,a}^{\alpha,\beta}$ by the weighted Korovkin type theorem given by Gadzhiev (\cite{Gad1}, \cite{Gad2}) and using weighted modulus of continuity defined by Ispir\cite{ispir}.\\ \\
Let $B_\rho[0,\infty)=\{f(x):|f(x)|\leq M_f \rho(x),\rho(x)$ is weight function, $M_f$  is a constant depending on $f$ and $x\in[0,\infty) \}$,\\
 $C_\rho[0,\infty)$ is the space of continuous function in $B_\rho[0,\infty)$ with the norm $\|f(x)\|_\rho=\sup\limits_{x\in[0,\infty)}\frac{|f(x)|}{\rho(x)}$ and \\
 $C_\rho^{k}=\{f\in C_\rho: \lim\limits_{|x|\rightarrow\infty}\frac{f(x)}{\rho(x)}=k,$ where $k$ is a constant depending on $f\}$.\\ \\
 \textbf{Theorem 4.1} Let $T_{n,a}^{\alpha,\beta}$ be the sequence of linear positive operators defined by (7). Then for $f\in C_\rho^k$,
 \begin{eqnarray*}
 \lim\limits_{n\rightarrow\infty}\|T_{n,a}^{\alpha,\beta}(f;x)-f(x)\|_\rho=0.
 \end{eqnarray*}
 \textbf{Proof} To prove the theorem, it is sufficient to show that \\
 \begin{eqnarray*}
 \lim\limits_{n\rightarrow\infty}\|T_{n,a}^{\alpha,\beta}(t^i;x)-x^i\|_\rho=0,\hspace{1cm}for\hspace{0.3cm} i=0,1,2.
 \end{eqnarray*}
 It is obvious that  $\lim\limits_{n\rightarrow\infty}\|T_{n,a}^{\alpha,\beta}(1;x)-1\|_\rho=0$. Now,from the Lemma 2.2, we have
 \begin{eqnarray*}
\|T_{n,a}^{\alpha,\beta}(t;x)-x\|&=&\sup\limits_{x\in[0,\infty)}\frac{\Big|\bigg(\frac{n}{n+\beta}-1\bigg)x+\frac{a}{n+\beta}\frac{x}{1+x}+ \frac{2\alpha+1}{2(n+\beta)}\Big|}{1+x^2}\\
 &&\leq\bigg(\frac{\beta}{n+\beta}\bigg)\sup_{x\in[0,\infty)}\frac{x}{1+x^2}+\frac{a}{n+\beta}\sup\limits_{x\in[0,\infty)}\frac{x}{(1+x)(1+x)^2}+\frac{2\alpha+1}{2(n+\beta)}\sup\limits_{x\in[0,\infty)}\frac{1}{(1+x)^2}.
 \end{eqnarray*}
 \begin{eqnarray*}
  \parallel \Hat{T}_n^{\alpha,\beta}(t;x)-x\parallel_\rho\rightarrow 0 \hspace{1cm}as\hspace{0.2cm} n\rightarrow\infty.
 \end{eqnarray*}
 Also, we can write
 \begin{eqnarray*}
 \| T_{n,a}^{\alpha,\beta}(t^2;x)-x^2\|&= &\sup\limits_{x\in[0,\infty)} \frac{\Big|\frac{n^2+n}{(n+\beta)^2}x^2+\frac{n(2+2\alpha)}{(n+\beta)^2}x+\frac{a^2}{(n+\beta)^2}\frac{x^2}{(1+x)^2}+\frac{2an}{(n+\beta)^2}\frac{x^2}{(1+x)}+\frac{a(2+2\alpha)}{(n+\beta)^2}\frac{x}{1+x} + \frac{3\alpha^2+1}{3(n+\beta)^2}-x^2\Big|}{1+x^2}\\
 &&\leq\frac{n(1-2\beta)+\beta^2}{(n+\beta)^2}\sup\limits_{x\in[0,\infty)}\frac{x^2}{1+x^2}+\frac{n(2+2\alpha)}{(n+\beta)^2}\sup\limits_{x\in[0,\infty)}\frac{x}{1+x^2}+\frac{a^2}{(n+\beta)^2}\sup_{x\in[0,\infty)}\frac{x^2}{(1+x)^2(1+x^2)}\\
 &&+\frac{2an}{(n+\beta)^2}\sup\limits_{x\in[0,\infty)}\frac{x^2}{(1+x)(1+x^2)}+\frac{a(2+2\alpha)}{(n+\beta)^2}\sup\limits_{x\in[0,\infty)}\frac{x}{(1+x)(1+x^2)}+ \frac{3\alpha^2+1}{3(n+\beta)^2}\sup\limits_{x\in[0,\infty)}\frac{1}{(1+x^2)}.
\end{eqnarray*}
Which implies that $\| T_{n,a}^{\alpha,\beta}(t^2;x)-x^2\|_\rho\rightarrow 0$ as $n\rightarrow\infty$.
Hence, proof is completed.\\ \\ \\
Now we prove the next theorem using weighted modulus of continuity defined by  Ispir\cite{ispir} as:
\begin{eqnarray*}
\Omega_n(f;\delta)=\sup\limits_{|h|\leq\delta,x\in [0,\infty)}\frac{|f(x+h)-f(x)|}{(1+h^2)(1+x^2)},
\end{eqnarray*}
for each $f\in C_\rho^k[0,\infty)$. \\
Some properties of $\Omega_n(f;\delta)$ are as follows:\\ \\
\textbf{Lemma 4.2}  Let $f\in C_\rho^K$,\\ \\
(i) $\Omega_n(f;\delta)$ is a monotonically increasing function of $\delta$, $\delta\geq 0$.\\ \\
(ii) For every $f\in C_\rho ^k, \lim\limits_{\delta\rightarrow0}\Omega_n(f;\delta)=0$.\\ \\
(iii) For each positive value of $\lambda$
\begin{eqnarray*}
 \Omega_n(f;\lambda\delta)\leq 2(1+\lambda)(1+\delta^2)\Omega_n(f;\delta).
\end{eqnarray*}
Using this property of modulus of continuity $\Omega_n(f;\delta)$ and its definition, we have
\begin{eqnarray*}
|f(t)-f(x)|\leq(1+x^2)(1+(t-x)^2)\Omega_n(f;|t-x|)
\end{eqnarray*}
and also
\begin{eqnarray}
 |f(t)-f(x)|&\leq2&\Big(\frac{|t-x|}{\delta_n}+1\Big)\Omega_n(f;\delta_n)(1+x^2)(1+(t-x)^2).\\
 \Big(\frac{|t-x|}{\delta_n}+1\Big)(1+(t-x)^2)&\leq2&(1+\delta_n^2)\Big(1+\frac{(t-x)^4}{\delta_n^4}\Big)
\end{eqnarray}
\textbf{Theorem 4.3} Let $f\in C_{\rho,[0,\infty)}^k$. Then the inequality
\begin{eqnarray*}
\sup\limits_{x\in [0,\infty)}\frac{|T_{n,a}^{\alpha,\beta}(f;x)-f(x)|}{(1+x^2)^3}\leq M_a^{\alpha,\beta}\Omega_n(f;(n+\beta)^{-\frac{1}{2}})
\end{eqnarray*}
holds, where $M_a^{\alpha,\beta}$ is a constant independent of n.\\ \\
\textbf{Proof} From (16) and (17) we have
\begin{eqnarray*}
|T_{n,a}^{\alpha,\beta}(f;x)-f(x)|&\leq& 4(1+x^2)(1+\delta_n^2)\Omega_n(f;\delta_n)T_{n,a}^{\alpha,\beta}\Big(1+\frac{(t-x)^4}{\delta_n^4}\Big),\\
 &=&4(1+x^2)(1+\delta_n^2)\Omega_n(f;\delta_n)\Big(1+T_{n,a}^{\alpha,\beta}\Big(\frac{(t-x)^4}{\delta_n^4};x\Big)\Big)\\
 &\leq& (1+x^2)\Omega_n(f;\delta_n)\Big(1+\frac{1}{\delta^4}T_{n,a}^{\alpha,\beta}((t-x)^4;x).\\
\end{eqnarray*}
Using Lemma 2.4, we obtain
\begin{eqnarray*}
|T_{n,a}^{\alpha,\beta}(f;x)-f(x)|&\leq& M_a^{\alpha,\beta}(1+x^2)\Omega_n(f;\delta_n)\Big(1+\frac{\Big(\frac{1}{(n+\beta^2)}\Big)(x^4+x^3+x^2+x+1)}{\delta_n^4}\Big)\\
 &\leq& M_a^{\alpha,\beta}(1+x^2)\Omega_n(f;(n+\beta)^{-\frac{1}{2}})(x^4+x^3+x^2+x+1)\\
 &\leq&M_a^{\alpha,\beta}(1+x^2)^3\Omega_n(f;(n+\beta)^{-\frac{1}{2}}).
\end{eqnarray*}
Therefore, we obtain
\begin{eqnarray*}
\sup\limits_{x\in [0,\infty)}\frac{|T_{n,a}^{\alpha,\beta}(f;x)-f(x)|}{(1+x^2)^3}\leq M_a^{\alpha,\beta}\Omega_n(f;(n+\beta)^{-\frac{1}{2}}).
 \end{eqnarray*}

\end{document}